\documentclass[reqno,12pt,a4paper]{amsart}

\usepackage{tikz}
\usetikzlibrary{matrix,arrows,calc,snakes,patterns,decorations.markings}

 \usepackage{graphicx}
\usepackage[mathscr]{eucal}
\usepackage{graphics,epic}
\usepackage{amsfonts}
\usepackage{amscd}
\usepackage{latexsym}
\usepackage{amsmath,amssymb, amsthm, stmaryrd, bm, bbm}
\usepackage[all,2cell]{xy}
\usepackage{mathrsfs}
\usepackage{color}

\newtheorem*{theorem*}{Theorem}

\newtheorem*{conjecture*}{Conjecture}

\newtheorem*{question*}{Question}
\theoremstyle{remark}

\theoremstyle{definition}

\newtheorem*{notation*}{Notation}
\newtheorem*{proof of theorem 1.2}{Proof of Theorem 1.2}

\numberwithin{equation}{section}

\begin{document}

\title[]{On mutation of $\tau$-tilting modules}

\author{Yingying Zhang}

\address{Department of Mathematics, Nanjing University, Nanjing 210093, jiangsu Province, P.R.China}

\email{zhangying1221@sina.cn}

\begin{abstract}
Mutation of $\tau$-tilting modules is a basic operation to construct a new support $\tau$-tilting module from a given one by replacing a direct summand. The aim of this paper is to give a positive answer to the question posed in [AIR, Question 2.31] about mutation of $\tau$-tilting modules.
\end{abstract}

\maketitle

\section{Introduction}\label{s:introduction}

$\tau$-tilting theory was introduced by Adachi, Iyama and Reiten [AIR] and completes (classical) tilting theory from the viewpoint of mutation. Note that $\tau$-tilting theory has stimulated several investigations; in particular, there is close relation between support $\tau$-tilting modules (see definition 2.1 for details) and some other important notions in representation theory, such as torsion classes, silting complexes, cluster-tilting objects and $*$-modules (see [AIR, AiI, AMV, IJY, IR, J, W] and so on). Since $\tau$-tilting theory was introduced, many algebraists started to apply it to important classes of algebras (see [Ad, EJR, HZ, IZ, M1, M2, Z] and so on).

Let us recall a main result in the paper [AIR]. Let $\Lambda$ be a finite dimensional algebra and $T=X\oplus U$ a basic $\tau$-tilting $\Lambda$-module with an indecomposable summand $X$ satisfying $X\notin$ ${ \rm Fac\,} U$. Take an exact sequence
\begin{equation}
X\buildrel {f} \over\longrightarrow U'\longrightarrow Y\longrightarrow 0
\end{equation}
with a left $({\rm add\,}U)$-approximation $f$. It is shown in [AIR, Theorem 2.30] that $Y$ is either zero or a direct sum of copies of an indecomposable $\Lambda$-module $Z$, and we can obtain a new basic support $\tau$-tilting $\Lambda$-module $\mu_{X}(T)$ called $mutation$ of $T$ with respect to X by $\mu_{X}(T)=U$ if $Y=0$ and $\mu_{X}(T)=Z\oplus U$ if $Y\neq 0$.

Moreover, they posed the following question.

\noindent{\bf Question 1.1.}
Assume that $Y$ in (1.1) is nonzero. Is $Y$ indecomposable?

\vspace{0.2cm}

A partial answer for the case when $\Lambda$ is an endomorphism algebra of a cluster-tilting object was given by Yang and Zhu in [YZ, Corollary 4.17].
The aim of this paper is to give a positive answer to this question.

\vspace{0.2cm}

\noindent{\bf Theorem 1.2.}
If $Y$ in (1.1) is nonzero, then it is indecomposable.

\vspace{0.2cm}

The idea of proof is to use the bijection between support $\tau$-tilting modules and two-term silting complexes given in [AIR].

\vspace{0.2cm}

\begin{notation*}
Let $K$ be a field. By an algebra $\Lambda$, we mean a finite dimensional $K$-algebra. We denote by $\mathrm{mod}\,\Lambda$ (resp. ${ \rm proj\,} $$\Lambda$) the category of finitely generated (resp. finitely generated projective) left $\Lambda$-modules and by $\tau$ the Auslander-Reiten translation of $\Lambda$. We denote by ${ \rm K}^{\rm b}({ \rm proj\,\Lambda})$ the homotopy category of bounded complexes of finitely generated projective $\Lambda$-modules. For $X\in\mathrm{mod}\,\Lambda$, we denote by ${ \rm add\,} X$ (resp. ${ \rm Fac\,} X$) the subcategory of $\mathrm{mod}\,\Lambda$ consisting of direct summands (resp. factor modules) of finite direct sums of copies of $X$.
\end{notation*}

\section{Proof of theorem}

First we recall the definition of support $\tau$-tilting modules from [AIR].

\vspace{0.2cm}

\noindent{\bf Definition 2.1.}
Let $X \in {\rm mod\,} \Lambda$ and $P \in {\rm proj\,} \Lambda$.
\begin{itemize}
\item[(1)] We call $X$ {\it
$\tau$-rigid} if ${\rm Hom}_{\Lambda}(X, \tau X) = 0$. We call ($X$,$P$) a $\tau$-$rigid$  $pair$ if $X$ is $\tau$-rigid and ${\rm Hom}_{\Lambda}(P, X)$=0.

\item[(2)] $X$ is called {\it $\tau$-tilting} if $X$ is
$\tau$-rigid and $|X| = |\Lambda|$.

\item[(3)] $X$ is called {\it support $\tau$-tilting} if there
exists an idempotent $e$ of $\Lambda$ such that $X$ is a $\tau$-tilting
$(\Lambda/\langle e\rangle)$-module. We call ($X$,$P$) a $support$ $\tau$-$tilting$ $ pair$ if ($X$,$P$) is $\tau$-rigid and $|X|+|P|=|\Lambda|$.
\end{itemize}

\vspace{0.2cm}
If ($X$,$P$) is a support $\tau$-tilting pair for $\Lambda$, then $X$ is a support $\tau$-tilting $\Lambda$-module. Conversely any support $\tau$-tilting $\Lambda$-module $X$ can be extended to a support $\tau$-tilting pair ($X$,$P$). We denote by {\rm s$\tau$-tilt$\,\Lambda$} the set of isomorphism classes of basic support $\tau$-tilting pairs (or equivalently, modules) for $\Lambda$.

\vspace{0.2cm}

Now we recall the definition of silting complexes from [AiI].

\vspace{0.2cm}

\noindent{\bf Definition 2.2.}
We call $P$ $\in$ ${ \rm K}^{\rm b}({ \rm proj\,\Lambda})$ $silting$ if ${\rm Hom}$$_{{ \rm K}^{\rm b}({ \rm proj\,\Lambda})}$($P$,$P$[$i$])=0 for any $i>0$ and ${\rm thick\,}$$P$= ${ \rm K}^{\rm b}({ \rm proj\,\Lambda})$, where ${\rm thick\,}$$P$ is the smallest full subcategory of ${ \rm K}^{\rm b}({ \rm proj\,\Lambda})$ containing $P$ and is closed under cones, [$\pm 1$] and direct summands. A complex $P$=($P^{i},d^{i}$) in ${ \rm K}^{\rm b}({ \rm proj\,\Lambda})$ is called $two\text{-}term$ if $P^{i}$=0 for all $i\neq0, -1$.

\vspace{0.2cm}

We denote by $2$-${\rm silt\,\Lambda}$ the set of isomorphism classes of basic two-term silting complexes in ${ \rm K}^{\rm b}({ \rm proj\,\Lambda})$. Moreover, recall the definition of mutation of silting complexes from [AiI].

\vspace{0.2cm}

\noindent{\bf Definition-Proposition 2.3.} [AiI, Theorem 2.31]
Let $P=X\oplus Q$ be a basic silting complex in ${ \rm K}^{\rm b}({ \rm proj\,\Lambda})$ with an indecomposable summand $X$. We take a minimal left $({\rm add\,}Q)$-approximation $f$ and a triangle
\[X\buildrel {f} \over\longrightarrow Q'\longrightarrow Y\longrightarrow X[1]. \] Then $Y$ is indecomposable and $\mu_{X}^{-}(P):=Y\oplus Q$ is a basic silting complex in ${ \rm K}^{\rm b}({ \rm proj\,\Lambda})$ called the {\it left mutation} of $P$ with respect to $X$. The {\it right mutation} is defined dually but we will not use it in this paper.
\vspace{0.2cm}

The following result establishes a relation between $2$-${\rm silt\,\Lambda}$ and {\rm s$\tau$-tilt$\,\Lambda$}.

\vspace{0.2cm}

\noindent{\bf Theorem 2.4.} [AIR, Theorem 3.2 and Corollary 3.9] There exists a bijection
\begin{equation}
2\text{-}{\rm silt\,\Lambda} \longleftrightarrow {\rm s\tau\text{-}tilt\,\Lambda}
\end{equation}
given by $2$-${\rm silt\,\Lambda}$ $\ni P\mapsto H^{0}(P) \in$ {\rm s$\tau$-tilt$\,\Lambda$} and {\rm s$\tau$-tilt$\,\Lambda$} $\ni (M,P)\mapsto (P_{1}\oplus P \buildrel {(f, 0)} \over \rightarrow P_{0}) \in$ $2$-${\rm silt\,\Lambda}$, where $f:P_{1}\rightarrow P_{0}$ is a minimal projective presentation of $M$.
Moreover, the bijection (2.1) preserves mutation.

\vspace{0.2cm}

Now we give the proof of main result of this paper.

\begin{proof}[Proof of Theorem 1.2.]
Assume $Y\neq 0$ in (1.1).

By taking the minimal projective presentations of $X$ and $U$. We have the following exact sequences: \[P_{X}^{-1}\buildrel {d_{X}^{-1}} \over\longrightarrow P_{X}^{0}\buildrel {d_{X}^{0}} \over\longrightarrow X\longrightarrow 0,\quad P_{U}^{-1}\buildrel {d_{U}^{-1}} \over\longrightarrow P_{U}^{0}\buildrel {d_{U}^{0}} \over\longrightarrow U\longrightarrow 0.\]
Denote by $P_{X}=(P_{X}^{-1}\buildrel {d_{X}^{-1}} \over\longrightarrow P_{X}^{0})$ and $P_{U}=(P_{U}^{-1}\buildrel {d_{U}^{-1}} \over\longrightarrow P_{U}^{0})$. Then $P_{T}=P_{X}\oplus P_{U}$ gives a minimal projective presentation of $T$. Then by Theorem 2.4 it follows that $P_{T}$ belongs to $2$-${\rm silt\,\Lambda}$ since $T$ is a basic $\tau$-tilting $\Lambda$-module. Also $P_{X}$ is indecomposable since $X$ is indecomposable.

Under above setting, $\mu_{X}(T)\in$ {\rm s$\tau$-tilt$\,\Lambda$} and $\mu_{P_{X}}^{-}(P_{T})\in$ $2$-${\rm silt\,\Lambda}$ correspond via the bijection (2.1) in Theorem 2.4. In particular, we have
\[\mu_{X}(T)=H^{0}(\mu_{P_{X}}^{-}(P_{T})).\]
To calculate $\mu_{P_{X}}^{-}(P_{T})$ in ${ \rm K}^{\rm b}({ \rm proj\,\Lambda})$, we take a triangle
\begin{equation}
P_{X}\buildrel {a} \over\longrightarrow P'\buildrel {b} \over\longrightarrow Q\longrightarrow P_{X}[1],
\end{equation}
where $a$ is a minimal left $({ \rm add\,} P_{U})$-approximation of $P_{X}$. Then by Definition-Proposition 2.3, $Q$ is indecomposable and $\mu_{P_{X}}^{-}(P_{T})=Q\oplus P_{U}$.

Taking the 0th cohomology of the triangle (2.2), we obtain the following exact sequence:
\[X\buildrel f'\over\longrightarrow U_{P'}\buildrel g' \over\longrightarrow Y_{Q}\longrightarrow 0,\]
where $U_{P'}=H^{0}(P')\in H^{0}({ \rm add\,} P_{U})={ \rm add\,} U$, $Y_{Q}=H^{0}(Q), f'=H^{0}(a)$ and $g'=H^{0}(b)$. Since $Q$ is indecomposable, it follows that $Y_{Q}$ is indecomposable.

We claim that $f'$ is a left $({ \rm add\,} U)$-approximation.
For any $h\in{ \rm Hom}_{\Lambda}(X, U)$, there exist morphisms $c^{-1}$ and $c^{0}$ making the following diagram commutative:
$$\xymatrix{
P_{X}^{-1}\ar[r]^{d_{X}^{-1}}\ar[d]^{c^{-1}}&P_{X}^{0}\ar[r]^{d_{X}^{0}}\ar[d]^{c^{0}}& X\ar[r]\ar[d]^{h}&0\\
P_{U}^{-1}\ar[r]_{d_{U}^{-1}}&P_{U}^{0}\ar[r]_{d_{U}^{0}}& U\ar[r]&0
.}
$$
Define $c=\overline{(c^{-1}, c^{0})}\in{ \rm Hom}_{{ \rm K}^{\rm b}({ \rm proj\,\Lambda})}(P_{X}, P_{U})$. Immediately, we have $H^{0}(c)=h$.

Since $a\in{ \rm Hom}_{{ \rm K}^{\rm b}({ \rm proj\,\Lambda})}(P_{X}, P')$ is a left $({ \rm add\,} P_{U})$-approximation, there exists $e\in{ \rm Hom}_{{ \rm K}^{\rm b}({ \rm proj\,\Lambda})}(P', P_{U})$ such that $c=ea\in{ \rm Hom}_{{ \rm K}^{\rm b}({ \rm proj\,\Lambda})}(P_{X}, P_{U})$. Since $H^{0}(-): { \rm K}^{\rm b}({ \rm proj\,\Lambda})\rightarrow{\rm mod\,}\Lambda$ is a functor, we have $$h=H^{0}(c)=H^{0}(ea)=H^{0}(e)H^{0}(a)=H^{0}(e)f'.$$ We have finished to prove that $f'$ is a left $({ \rm add\,} U)$-approximation.

Since $f$ is a minimal left $({ \rm add\,} U)$-approximation, there exists a module $U''$ in ${ \rm add\,} U$ such that $U_{P'}=U'\oplus U''$ and $Y_{Q}=Y\oplus U''$. Since $Y_{Q}$ is indecomposable and $Y\neq 0$ by our assumption, we have that $U''=0$ and $Y=Y_{Q}$ is indecomposable.
\end{proof}

\vspace{0.5cm}

{\bf Acknowledgement}
This paper was written while the author was visiting Nagoya University from October 2015 to September 2016. The author thanks Prof. Osamu Iyama for the careful guidance and Prof. Zhaoyong Huang for the continuous encouragement. She thanks Yuya Mizuno for his valuable comments. She also wants to thank people in Nagoya University for their help. This work was partially supported by NSFC (Grant No. 11571164) and the China Scholarship Council (CSC).

\end{document}